\documentclass[a4paper]{jpconf}
\usepackage{graphicx}

\usepackage{amssymb}
\usepackage{amsmath}
\usepackage{amsthm}
\usepackage{amsfonts}

\newtheorem{theorem}{Theorem}[section]
\newtheorem{corollary}[theorem]{Corollary}

\newtheorem{example}[theorem]{Example}

\newtheorem{remark}[theorem]{Remark}

\newtheorem{definition}[theorem]{Definition}
\numberwithin{equation}{section}

\def\O{\Omega}

\def\s{\sigma}

\def\l{\lambda}

\def\P{\Phi}
\def\O{\Omega}

\def\L{\Lambda}

\begin{document}
\title{On cubic stochastic operators and processes
\footnote{To the memory of our teacher T. A. Sarymsakov on the occasion of his 100th birthday.}}

\author{B. J. Mamurov$^1$, U. A. Rozikov$^2$}

\address{$^1$Bukhara State Medical Institute,
 1, Navoi str., \, Bukhara, Uzbekistan\\
 $^2$ Institute of Mathematics, National University of Uzbekistan, 29, Do'rmon Yo'li str., 100125,
Tashkent, Uzbekistan} \ead{bmamurov.51@mail.ru,
rozikovu@yandex.ru}

\begin{abstract}  In this paper analogically as quadratic stochastic operators and processes we
define cubic stochastic operator (CSO) and cubic stochastic processes (CSP).
These are defined on the set of all probability
measures of a measurable space. The measurable space can be given on a finite or continual set.
The finite case has been investigated before. So here we mainly work on the continual set.
We give a construction of a CSO and show that dynamical systems generated by such a CSO can be
studied by studying of the behavior of trajectories of a CSO given on a finite dimensional simplex.
We define a CSP and drive differential equations for such CSPs with continuous time.
\end{abstract}

%\maketitle

%{\bf Mathematics Subject Classifications (2010).} 37C20; 37C25.

%{\bf{Key words.}} Cubic stochastic operator; cubic process; trajectory.

\section{Introduction}

Dynamical systems generated by quadratic operators have been proved to be a rich
source of analysis for the study of dynamical properties and
modeling in different domains, such as population dynamics
\cite{B,FG}, physics \cite{UR}, economy \cite{D}, mathematics
\cite{GR1},\cite{GN},\cite{ES},\cite{J},\cite{ly}. The approach to population genetics, posed
within that scheme the problem of an explicit description of
evolutionary operators of free populations.

A {\it quadratic} stochastic operator (QSO) of a free population \cite{ly} is a
(quadratic) mapping of the simplex
$$ S^{m-1}=\{x=(x_1,...,x_m)\in R^m: x_i\geq 0, \sum^m_{i=1}x_i=1 \} $$
into itself, of the form
\begin{equation}\label{1}
 V: x_k'=\sum^m_{i,j=1}p_{ij,k}x_ix_j, \ \ (k=1,...,m)
\end{equation}
where $p_{ij,k}$ are coefficients of heredity and
\begin{equation}\label{2}
 p_{ij,k}\geq 0, \ \ \sum^m_{k=1}p_{ij,k}=1, \ \ (i,j,k=1,...,m)
\end{equation}
 Note that each element $x\in S^{m-1}$ is a probability
distribution on $E=\{1,...,m\}.$

Similarly, one can define a {\it cubic} stochastic operator (CSO) $W:S^{m-1}\to S^{m-1}$ as
\begin{equation}\label{c1}
 W: x_l'=\sum^m_{i,j,k=1}P_{ijk,l}x_ix_jx_k, \ \ (l=1,...,m)
\end{equation}
where $P_{ijk,l}$ are coefficients such that
\begin{equation}\label{c2}
 P_{ijk,l}\geq 0, \ \ \sum^m_{l=1}P_{ijk,l}=1, \ \ (i,j,k,l=1,...,m)
\end{equation}

The population evolves by starting from an arbitrary state (probability
distribution on $E$) $x\in S^{m-1}$ then passing to the state $Vx$ (or $Wx$)
(in the next ``generation''), then to the state $V^2x$ (resp. $W^2x$), and so on.
 Thus, states of the population described by the following discrete-time dynamical
 system
 \begin{equation}\label{2}
 x^{0},\ \ x'= V(x), \ \ x''=V^{2}(x),\ \  x'''= V^{3}(x),\dots
\end{equation}
where $V^n(x)=\underbrace{V(V(...V}_n(x))...)$ denotes the $n$ times
iteration of $V$ to $x$.
The main problem for a given dynamical system (\ref{2}) is to
describe the limit points of $\{x^{(n)}\}_{n=0}^\infty$ for
arbitrary given $x^{(0)}$.

In \cite{GMR} (see also \cite{MG}) a review of the theory of QSOs
is given. Note that each quadratic (resp. cubic) stochastic
operator can be uniquely defined by a stochastic matrix ${\mathbf
P}= \{p_{ij,k}\}^m_{i,j,k=1}$ (resp. ${\mathbb
P}=\{P_{ijk,l}\}^m_{i,j,k,l=1}$). In \cite{GaR} a constructive
description of ${\mathbf P}$ (i.e. a QSO) is given. This
construction depends on cardinality of $E$, which can be finite or
continual. Some particular cases of this construction were defined
in \cite{GN}. In \cite{RK} a similar construction for the CSOs on
a finite set $E$ is given. We note that for the continual set $E$
one of the key problem is to determine the set of coefficients of
heredity which is already infinite dimensional. In this paper we
shall give a construction of CSO for a continual set $E$. In this
construction the CSO depends on a probability measure $\mu$ given
on a measurable space $(E,{\mathcal F})$. We will show that the
dynamical systems generated by the CSO depending on $\mu$ can be
reduced to a dynamical system generated by a Volterra CSO defined
on a finite-dimensional simplex.

Next goal of the paper is to define cubic stochastic processes
(CSP) and give differential equations for such processes.  This
investigations will be similar to works
\cite{G},\cite{MS2015,MS2013}, \cite{S1}-\cite{S3}, where the
authors introduced a continuous-time dynamical system as quadratic
stochastic processes (QSP). The reader is referred to very recent
book \cite{MG} for the theory of QSPs.

\section{Definitions and examples}

Consider a measurable space $(E, {\mathcal F})$ and let $S(E,{\mathcal F})$
be the set of all probability measures on $(E,{\mathcal F}).$

\begin{definition}\label{d1}
A mapping $W : S(E,{\mathcal F})\to S(E,{\mathcal F})$ is called a
cubic stochastic operator (CSO) if, for an arbitrary measure
$\l\in S(E,{\mathcal F}),$ the measure $\l'=W\l$ is defined by
$$ \l'(A)=\int_E\int_E\int_E P(x, y, z, A)d\l(x)d\l(y)d\l(z), \ \ \forall A\in {\mathcal F},$$
where $P(x, y, z, A)$, satisfies the following conditions:

(i) $P(x, y, z,\cdot)\in S(E,{\mathcal F})$ for any fixed $x,y,z \in E;$

(ii) $P(x, y, z, A),$ regarded a function of three variables $x$, $y$, and $z$
is measurable on $(E\times E\times E, {\mathcal F}\otimes {\mathcal F}\otimes {\mathcal F})$ for any fixed
$A\in {\mathcal F}.$
\end{definition}

When $E$ is finite, a CSO on $S(E,{\mathcal F})=S^{m-1}$ is defined as
in (\ref{c1}) with $P_{ijk,l}=P(i,j,k,l).$

Let $ \left( {E, \mathcal F, {\mathcal M}} \right)$ be triple, where $\mathcal F$ is a $ \sigma $
-algebra of subsets of $E$ and ${\mathcal M}=S(E,{\mathcal F})$, i.e. the set of all probability measures on
$(E, \mathcal F)$. For any three elements $x, y, z\in E $ and a given measure $m_{t_0 } \in {\mathcal M}$
 at some moment $t_0$ of the time we assume that we know the law of probability distribution $ m_{t_1} \in {\mathcal M}$ of the system
 $E$ at the moment $ t_1 > t_0 $ of time.

 Denote by $ P\left( {t_1, x, y, z, t_2, A} \right)$ the probability of obtaining an element from the set $ A \in \mathcal F$ at the moment $t_2$,
 provided that the elements
$ x, y $ and $ z $ of $E$ interact starting at moment $t_1$, where $t_2  \ge t_1  + 1$.

Thus if at the moment $t_1$ we start with a probability measure $m_{t_1}\in {\mathcal M}$, then
$m_{t_2}\in {\mathcal M}$ for any $t_2\geq t_1+1$ is defined by
\begin{equation}\label{m}
m_{t_2}(A)=\int_E\int_E\int_E P(t_1, x, y, z, t_2, A)m_{t_1}(dx)m_{t_1}(dy)m_{t_1}(dz).
\end{equation}
Without loss of generality we assume that the process starts at the moment $t=0$.

\begin{definition}\label{d2}
A family $ \{P\left(t_1, x, y, z, t_2, A\right)$: $x, y, z \in E$,
$A\in \mathcal F$, $t_1, t_2\in \mathbb R^+$, $t_2-t_1\geq 1\}$ is
called {\rm cubic stochastic process} (CSP) if it satisfies the
following conditions
\begin{itemize}
\item[(I)] $ P\left(t, x, y, z, t + 1, A\right) = P\left(0, x, y, z, 1, A\right)$ for any $ t \ge 1 $;

\item[(II)] The value of $P\left(t_1, x, y, z, t_2, A\right)$ is independent on any permutations of
variables $x, y, z$ for all $x, y, z \in E $, and  $ A \in \mathcal F $;

\item[(III)]  $ P\left(t_1, x, y, z, t_2, A\right)$ is a probability measure on $(E, \mathcal F)$ for all $ x, y, z \in E$ and $t_1, t_2\in \mathbb R^+$, $t_2-t_1\geq 1$;

\item[(IV)] $ P\left(t_1, x, y, z, t_2, A\right)$ as function of the three variables
$x, y$ and $z$ is measurable with respect to $\left(E\times E \times E,\quad \mathcal F \otimes \mathcal F \otimes \mathcal F\right)$ for all
$ A \in \mathcal F $;

\item[(V)] For any $ t_1  < t_2  <t_3 $ such that $ t_2  - t_1
\ge 1 $ and $ t_3  - t_2  \ge 1 $ the following holds
\begin{equation}\label{b1}
P\left(t_1, x, y, z, t_3, A\right) = \int\limits_E \int\limits_E
\int\limits_E P\left(t_1, x, y, z, t_2, du\right)P\left(t_2, u, \vartheta, q, t_3, A\right)m_{t_2}(d\vartheta)m_{t_2}(dq),
\end{equation}
where $m_{t_2}$ on $(E, \mathcal F)$ is defined by (\ref{m}).
\end{itemize}
\end{definition}

\begin{remark} Now we point out the followings:
\begin{itemize}
\item[1.] The consideration of CSO and CSP are motivated by their
appearance in biology (for example in gene engineering, a triple
crossings for different sorts of plants to obtain another sort,
and free population with à ternary production) \cite{A}, in
physics (for example in spin systems with ternary interactions)
\cite{mal}. We note that one may have more general form of
stochastic operators (resp. processes) with order $n\geq 1$, and
coefficients $P_{i_1i_2\dots i_n,j}\geq 0$ (resp. $P\left(t_1,
x_1, x_2,\dots, x_n, t_2, A\right)$) in biology (gene engineering)
and physics ($n$-nary interactions). In this paper, we restrict
ourselves to the case $n=3$, i.e. CSOs and CSPs. Because even for
the considered setting associated dynamical systems are very
complicated in comparison with $n=1$ (i.e. Markov chains) and
$n=2$ (i.e. QSOs).

\item[2.] Note that if one defines a new process by
$$
Q(s,x,y,t,A)=\int_E P(s,x,y,z,t,A)m_s(dz) $$ then one can check
that the defined process is   QSP (see \cite{MS2015,S1}). We
recall that the functions $Q(s,x,y,t,A)$ denote the probability
that under the interaction of the elements $x$ and $y$ at time $s$
an event $A$ comes into effect at time $t$. Since for physical,
chemical and biological phenomena, a certain time is necessary for
the realization of an interaction, it is taken the greatest such
time to be equal to 1 (see the Boltzmann model \cite{J} or the
biological model \cite{ly}). Thus the probability $Q(s,x,y,t,A)$
is defined for $ t-s \geq\ 1.$ Hence, for CSP, the probabilities
$P(s,x,y,z,t,A)$ are also defined for $ t-s \geq\ 1$.

\item[3.] It should be noted that the CSPs are related to CSOs
(see \cite{J,MG}) in the same way as Markov processes are related
to linear transformations.

\item[4.] The equation (\ref{b1}) in Definition \ref{d2} is an
analogue of Chapman-Kolmogorov equation. In \cite{G},
\cite{MS2013},\cite{S1}-\cite{S3} such equations were extended to
QSPs. We note that for QSPs there are two types of the
Chapman-Kolmogorov equations: type $A$ and type $B$. Similarly,
one also can define (at least) two types of the Chapman-Kolmogorov
equations for a CSP: The equation (\ref{b1}) corresponds to the
type $A$ of QSP, the following equation is an analogue of the type
$B$ for CSP:
\begin{eqnarray*}
P\left(s, x, y, z, t, A\right) &=& \int\limits_E \int\limits_E
\int\limits_E \int\limits_E \int\limits_E \int\limits_E P\left(s,
x, y_1, z_1, \tau, du\right)P\left(s, y, y_2, z_2, \tau,
dv\right)\\[2mm]
&&\times P\left(s, z, y_3, z_3, \tau, dw\right) P\left(\tau, u, v,
w, t, A\right)m_{s}(dy_1)
m_{s}(dz_1)\\[2mm]
&&m_{s}(dy_2)m_{s}(dz_2)m_{s}(dy_3)m_{s}(dz_3).
\end{eqnarray*}
 In this paper, for the sake of simplicity, we shall only consider CSPs which satisfy (\ref{b1}).
 \end{itemize}
 \end{remark}

Let us provide some examples of CSPs.

 Let $E = \left\{ {1,2,...,n} \right\} $ and $ x^{\left( 0 \right)}  = \left\{ {x_1 ^{\left( 0 \right)} ,x_2
^{\left( 0 \right)} ,...,x_n ^{\left( 0 \right)} } \right\} $ be an initial distribution on $E$.

Denote
$$ P_{ijk,l}:= P\left( {0,i,j,k,1,\{l\}} \right) ,\quad P_{ijk,l}^{\left[ {s,t} \right]}:=P\left(
{s,i,j,k,t,\{l\}} \right) . $$

By equation (\ref{m}) at the moment $ t = 1$ the vector $ x^{(1)}
= \left( {x_1^{(1)}, x_2^{(1)}, ... , x_n^{(1)} }\right)$ is
defined as follows
$$ \quad x_l^{\left( 1 \right)}  =
\sum\limits_{i,j,k = 1}^n {P_{ijk,l} x_i^{\left( 0 \right)}
x_j^{\left( 0 \right)} x_k^{\left( 0 \right)} } ,\quad l =
\overline {1,n}.$$ In this case the condition (I) reduces to $
P_{ijk,l}^{\left[ {t,t + 1} \right]} = P_{ijk,l}.$

In general, from (\ref{m}),(\ref{b1}) one finds
\begin{equation}\label{bu3} x_l^{(t)}  =
\sum\limits_{i,j,k} {P_{ijk,l}^{\left[ {s,t} \right]} x_i^{(s)},
x_j^{(s)} x_k^{(s)} }
\end{equation}
\begin{equation}\label{bu4} P_{ijk,l}^{\left[ {s,t} \right]} = \sum\limits_{m,\gamma ,\delta
} {P_{ijk,m}^{\left[ {s,\tau } \right]} P_{m\gamma \delta
,l}^{\left[ {\tau ,t} \right]} x_\gamma ^{(\tau )} x_\delta
^{(\tau )}},
\end{equation} where $ \tau  - s \ge 1 $ and $ t-\tau \ge
1.$
\begin{example} To define a CSP for $E=\{1,2\}$ we first define a CSO by the matrix $\left(P_{ijk,l}\right)$, where
$$P_{111,1}=1, \ \ P_{112,1}=P_{121,1}=P_{211,1}={2\over 3},\ \ P_{122,1}=P_{212,1}=P_{221,1}={1\over 3}, \ \ P_{222,1}=0;$$
$$P_{111,2}=0, \ \ P_{112,2}=P_{121,2}=P_{211,2}={1\over 3},\ \ P_{122,2}=P_{212,2}=P_{221,2}={2\over 3}, \ \ P_{222,2}=1.$$
It is easy to see that this CSO is the identity mapping, i.e.
$x_1'=x_1$, $x_2'=x_2$. Take an initial vector $(x,1-x)$ and using
this CSO by formulas (\ref{bu3}),(\ref{bu4}) one can define a CSP:
$$\begin{array}{lllll}
P_{111,1}^{[s,t]}={1\over 3^{t-s-1}}\left[1+(3^{t-s-1}-1)x\right],\\[3mm]
P_{112,1}^{[s,t]}=P_{121,1}^{[s,t]}=P_{211,1}^{[s,t]}={1\over 3^{t-s-1}}\left[{2\over 3}+(3^{t-s-1}-1)x\right],\\[3mm]
P_{122,1}^{[s,t]}=P_{212,1}^{[s,t]}=P_{221,1}^{[s,t]}={1\over 3^{t-s-1}}\left[{2\over 3}+(3^{t-s-1}-1)(1-x)\right],\\[3mm]
P_{222,1}^{[s,t]}={3^{t-s-1}-1\over 3^{t-s-1}}x,\\[3mm]
P_{ijk,2}^{[s,t]}=1-P_{ijk,1}^{[s,t]}, \ \ \mbox{for all} \ \ i,j,k\in E=\{1,2\}.
\end{array}$$
\end{example}

\begin{example} Let $E=\{1,2,\dots,n\}$. Take a family of stochastic vectors: $a(t)=(a_1(t), a_2(t), \dots, a_n(t))$,
i.e. $a_i(t)\geq 0$, $\sum_i a_i(t)=1$ for any $t\geq 0$. For each
pari $s,t$ we define a stochastic matrix
$Q^{[s,t]}=(q_{il}^{[s,t]})_{i,l\in E},$ where
$q_{il}^{[s,t]}=a_l(t)$ for all $i\in E$, i.e. it does not depend
on $s$.  It is easy to see that this matrix satisfies the
Kolmogorov-Chapman equation:
$$Q^{[s,t]}=Q^{[s,\tau]}Q^{[\tau,t]}, \ \ \mbox{for all} \ \ 0\leq s<\tau<t.$$
Now define functions
$$P(s,i,j,k,t,\{l\})=q_{il}^{[s,t]}=a_l(t).$$
One can check that the defined family $\{P(s,i,j,k,t,\{l\})\}$ is
a CSP.
\end{example}
\begin{example} (cf. Example 4.2.1 of \cite{MG}) Let $(E,{\mathcal F})$ be a measurable space and $m_0$ be an initial measure on this space. Consider the following
functions
$$P(s,x,y,z,t,A)={1\over 3^{t-s-1}}\left({\delta_x(A)+\delta_y(A)+\delta_z(A)\over 3}+(3^{t-s-1}-1)m_0(A)\right),$$
where $t-s\geq 1$, $x,y,z\in E$ and $A\in \mathcal F$, $$\delta_x(A)=\left\{\begin{array}{ll}
1, \ \ \mbox{if} \ \ x\in A\\
0, \ \ \mbox{if} \ \ x\notin A.
\end{array}
\right.$$ It is easy to see that the defined family is CSP.
\end{example}

\section{A construction of CSOs}

Recall that a construction of CSO for finite $E$ was given in
\cite{RK}. Therefore, in this section we consider the continual
case.

Let $G=(\L,L)$ be a countable graph. For a finite set $\Phi$
denote by $\O$ the set of all functions $\s:\L\to\P.$ Let
$S(\O,\P)$ be the set  of all probability measures defined on
$(\O, {\mathcal F}),$ where ${\mathcal F}$ is the $\s-$algebra
generated by the finite-dimensional cylindrical set. Let $\mu$ be a measure
on $(\O, {\mathcal F})$ such that $\mu(B)>0$ for any
finite-dimensional cylindrical set $B\in {\mathcal F}.$

Let $M\subset \L$ be a finite connected subgraph.
Two elements $\s, \varphi\in \O$ are called equivalent if
$\s(x)=\varphi(x)$ for any $x\in M$, i.e.  $\s(M)=\varphi(M).$ Let
$\xi=\{\O_i, i=1,2,...,|\P|^{|M|}\},$ be the partition of $\O$
generated by this equivalent relation, where $|\cdot |$ denotes
the cardinality of a set and $\O_i$ contains all equivalent
elements\footnote{We note that $\xi$ depends on $M$, therefore all
quantities which we define using $\xi$ also depend on $M$.
But for simplicity of formulas we will omit $M$ from our formulas.}.

Denote $\langle ijk,l \rangle=\delta_{il}+\delta_{jl}+\delta_{kl}$, where $\delta$
is the Kronecker's symbol, i.e.
$$\delta_{ij}=\left\{\begin{array}{ll}
1, \ \ \mbox{if} \ \ i=j\\
0, \ \ \mbox{if} \ \ i\ne j.
\end{array}
\right.
$$

Consider
\begin{equation}\label{5}
 P_{\s_1\s_2\s_3, \s}=\frac{\langle ijk,l \rangle \mu(\Omega_l)}{\mu(\Omega_i)+\mu(\Omega_j)+\mu(\Omega_k)} \ \ {\rm if} \ \ \s_1\in \O_i, \, \s_2\in
\O_j, \, \s_3\in \O_k, \, \s\in \O_l.
\end{equation}

For arbitrary $\s$ it is easy to see that $P_{\s_1\s_2\s_3, \s}$ is invariant with respect to any permutations of $\s_1, \s_2, \s_3$.

Then the coefficients $P(\s_1, \s_2, \sigma_3, A)$ \ \ $(\s_1, \s_2, \s_3\in \O, \ \ A\in
{\mathcal F})$ are defined as
$$
P(\s_1,\s_2, \s_3, A)=Z(\s_1,\s_2, \s_3)\int_A P_{\s_1\s_2\s_3,\s}d\mu(\s)=
Z(\s_1,\s_2, \s_3)\sum_{l=1}^m\int_{A\cap \Omega_l}P_{\s_1\s_2\s_3,\s}d\mu(\s),$$
where $m=|\Phi|^{|M|}$ and $Z(\s_1, \s_2, \s_3)$ is the normalizing factor, which is chosen by the condition
that $ P(\s_1,\s_2, \s_3, \Omega)=1$.

It is easy to obtain the following
\begin{equation}\label{6}
P(\s_1,\s_2, \s_3, A)=\left\{\begin{array}{lll}
\frac{\mu(\O_i)\mu(A\cap\O_i)+\mu(\O_j)\mu(A\cap \O_j)+\mu(\O_k)\mu(A\cap\O_k)}
{\mu^2(\Omega_i)+\mu^2(\Omega_j)+\mu^2(\Omega_k)} ,\\[2mm]
\ \ \ \ \ \ \ \ \ \ {\rm if} \ \
\s_1\in \O_i, \s_2\in \O_j, \s_3\in \O_k, i\ne k, j\ne k, i\ne j\\[3mm]
\frac{2\mu(\O_i)\mu(A\cap\O_i)+\mu(\O_j)\mu(A\cap \O_j)}
{2\mu^2(\Omega_i)+\mu^2(\Omega_j)}, \ \ {\rm if} \ \
\s_1,\s_2\in \O_i, \s_3\in \O_j, i\ne j\\[3mm]
\frac{\mu(A\cap \O_i)}{\mu(\Omega_i)}  \ \ {\rm if} \ \ \s_1, \s_2, \s_3\in
\O_i.\end{array}\right. \hfill
\end{equation}

The CSO $W$ acting on the set $S(\O, \P)$ is determined by
coefficients (\ref{6}) is defined as follows: for an arbitrary
measure $\l\in S(\O, \P),$ the measure $\l'=W\l$ is
\begin{equation}\label{8}
\l'(A)=\int_\O\int_\O\int_\O P(\s_1,\s_2, \s_3, A)d\l(\s_1)d\l(\s_2)d\l(\s_3).
\end{equation}
Using (\ref{6}) from (\ref{8}) we obtain
$$
\l'(A)=\sum_{i=1}^{m}a_i(A)\l^3(\O_i)+$$
\begin{equation}\label{cso}
3\sum_{i=1}^m\sum_{{j=1\atop j\ne i}}^mb_{ij}(A)\l^2(\O_i)\l(\O_j)+6\sum_{1\leq i<j<k\leq m}c_{ijk}(A)\l(\O_i)\l(\O_j)\lambda(\O_k),
\end{equation}
where
\begin{equation}\label{abc}
\begin{array}{lll}
a_i(A)={\mu(A\cap \O_i)\over \mu(\O_i)},\\[3mm]
b_{ij}(A)=\frac{2\mu(\O_i)\mu(A\cap\O_i)+\mu(\O_j)\mu(A\cap \O_j)}
{2\mu^2(\Omega_i)+\mu^2(\Omega_j)},\\[3mm]
c_{ijk}(A)=\frac{\mu(\O_i)\mu(A\cap\O_i)+\mu(\O_j)\mu(A\cap \O_j)+\mu(\O_k)\mu(A\cap\O_k)}
{\mu^2(\Omega_i)+\mu^2(\Omega_j)+\mu^2(\Omega_k)}.
\end{array}
\end{equation}
It is easy to see that
\begin{equation}\label{ao}
\begin{array}{lll}
a_i(\O_l)=\left\{\begin{array}{ll}
1, \ \ \mbox{if} \ \ l=i,\\
0, \ \ \mbox{if} \ \ l\ne i
\end{array}\right. \ \ \
b_{ij}(\O_l)=\left\{\begin{array}{lll}
\frac{2\mu^2(\O_i)}
{2\mu^2(\Omega_i)+\mu^2(\Omega_j)}, \ \ l=i\\[3mm]
\frac{\mu^2(\O_j)}
{2\mu^2(\Omega_i)+\mu^2(\Omega_j)}, \ \ l=j\\[3mm]
0, \ \ \ \ \ l\ne i, \, l\ne j.
\end{array}\right.\\[3mm]
c_{ijk}(\O_l)=\left\{\begin{array}{ll}
\frac{\mu^2(\O_l)}{\mu^2(\Omega_i)+\mu^2(\Omega_j)+\mu^2(\Omega_k)}, \ \ l\in\{i,j,k\}\\[3mm]
0, \ \ \ \ \ l\notin \{i,j,k\}.
\end{array}\right.
\end{array}
\end{equation}

For a given measure $\l\in S(\L,\P)$ the trajectory $\{\l^{(n)}\},
n=1,2,...$ of the
operator ~(\ref{8}) is defined by $\l^{(n+1)}(A)=W(\l^{(n)})(A),$
where $n=0,1,2,... $ and $\l^{(0)}=\l,$ $A\in {\mathcal F}.$

By ~(\ref{cso}) and (\ref{ao}) we have
$$
\l'(\O_l)=\sum_{i=1}^{m}a_i(\O_l)\l^3(\O_i)+$$
\begin{equation}\label{cl}
3\sum_{i=1}^m\sum_{{j=1\atop j\ne i}}^mb_{ij}(\O_l)\l^2(\O_i)\l(\O_j)+6\sum_{1\leq i<j<k\leq m}c_{ijk}(\O_l)\l(\O_i)\l(\O_j)\lambda(\O_k).
\end{equation}
A CSO (\ref{c1}) is called
Volterra if the coefficients $P_{ijk,l}$ may be nonzero only when
$l\in\{i,j,k\}$ and vanish in all the remaining cases (see \cite{Kh}, \cite{Kh1}).

It is easy to see that any Volterra CSO has the following form
\begin{equation}\label{vc}
W: \l_l'=\l_l\left(\l_l^2+\l_l\sum^m_{i=1\atop i\ne l}a_{i,l}\l_i+\sum^m_{i,j=1\atop i\ne l, \, j\ne l}b_{ij,l}\l_i\l_j\right), \ \ (l=1,...,m)
\end{equation}
where $a_{i,l}$ and $b_{ij,l}$ are some coefficients depending on $P_{ijk,l}$.

Denoting $\l_i=\l(\O_i)$, and $a_{i,l}=3b_{li}(\O_l)$, $b_{ij,l}=6c_{ijl}(\O_l)$ the operator ~(\ref{cl}) can be written as (\ref{vc}).

Note that the $n-$th iteration $\l^{(n)}=W^{(n)}\l^{(0)}$ of
the operator  ~(\ref{8}) (i.e. (\ref{cso})) can be written as
$$
\l^{(n+1)}(A)=\sum_{i=1}^{m}a_i(A)(\l_i^{(n)})^3+$$
\begin{equation}\label{cson}
3\sum_{i=1}^m\sum_{{j=1\atop j\ne i}}^mb_{ij}(A)(\l_i^{(n)})^2\l_j^{(n)}+6\sum_{1\leq i<j<k\leq m}c_{ijk}(A)\l^{(n)}_i\l^{(n)}_j\lambda^{(n)}_k,
\end{equation}
where $\l^{(n)}_j,$ $j=1,...,m$ are coordinates of the trajectory
of operator ~(\ref{vc}) for the given $\l.$

Thus in order to study the trajectory of operator~(\ref{8}) it is
enough to know the behavior of trajectories of the operator
~(\ref{vc}), i.e.
we proved the following

\begin{theorem}\label{t1} For any finite $M\subset \Lambda$ the dynamical system generated by the CSO (\ref{8})
is reducible to a dynamical system generated by a Volterra CSO acting on $(m-1)$-dimensional simplex.
\end{theorem}

From Theorem \ref{t1} we get
\begin{corollary}  Assume for a given measure $\mu$ and $\l=(\l_1,...,\l_m)\in S^{m-1}$ for the trajectory of the Volterra operator (\ref{vc}) we have
$$
\lim_{n\to\infty}\l^{(n)}=(\l_1^*,\l^*_2,...,\l^*_m).$$
Then the corresponding trajectory
$\{\l^{(n)}(A)\}$ of the operator ~(\ref{8}) has the following limit
$$
\l(A)=\lim_{n\to\infty}
\l^{(n)}(A)=\sum_{i=1}^{m}a_i(A)(\l_i^*)^3+$$
\begin{equation}\label{16}
3\sum_{i=1}^m\sum_{{j=1\atop j\ne i}}^mb_{ij}(A)(\l_i^*)^2\l_j^{*}+6\sum_{1\leq i<j<k\leq m}c_{ijk}(A)\l^*_i\l^*_j\lambda^*_k.
\end{equation}
\end{corollary}

\begin{remark} As it was mentioned above the theory of Volterra QSOs is well studied (see for example \cite{GMR}).
But Volterra CSOs  were not exhaustively
studied, because such cubic operators are still complicated.
There are just a few articles devoted to Volterra CSOs \cite{Kh}, \cite{Kh1}, \cite{RK1}, \cite{RK}.
Therefore formula (\ref{16}) is already helpful by using the results for the Volterra CSOs studied in these papers.
 \end{remark}

\section {Integro-differential equations for CSP}

The equations which we want to drive here were given in \cite{M} for finite $E$. So we consider the continual case of $E$.
Let $E$ be a continual set. Consider a CSP on a measurable space $ (E, \mathcal F) $ with initial measure
$ m_0 $. For $ t > s + 2 $ from condition (V) of Definition \ref{d2} we get
$$P(s,x,y,z,t + \Delta ,A) - P(s,x,y,z,t,A)$$
$$
= \int\limits_E \int\limits_E \int\limits_E P(s,x,y,z,t -
1,du)\left\{ P(t - 1,u,\vartheta ,q,t + \Delta ,A)\right.$$ $$\left. \ \ \ \ - P(t -
1,u,\vartheta ,q,t,A)\right\}m_{t - 1} (d\vartheta )  m_{t -
1} (dq).$$
Assume the following limit exists
$$ C(t,u,\vartheta ,q,A) = \mathop {\lim
}\limits_{\Delta  \to 0} \frac{{P(t - 1,u,\vartheta ,q,t + \Delta
,A) - P(t - 1,u,\vartheta ,q,t,A)}}{\Delta }.$$ Then taking limit $ \Delta  \to 0 $
we obtain the {\it first integro-differential equation}:
\begin{equation}\label{b10}
\frac{{\partial P(s,x,y,z,t,A)}}{{\partial t}} = \int\limits_E
{\int\limits_E {\int\limits_E {P(s,x,y,z,t - 1, du)\;C(t,u,\vartheta
,q,A)m_{t - 1} (d\vartheta )} } } m_{t - 1} (dq).
\end{equation}
Similarly, one gets the {\it second integro-differential equation}:
\begin{equation}\label{b11}
\frac{{\partial P(s,x,y,z,t,A)}}{{\partial s}} =  -
\int\limits_E {\int\limits_E {\int\limits_E {C(s + 1,x,y,z,du)P(s
+ 1,u,\vartheta ,q,t,A)m_{s + 1} (d\vartheta )} } } m_{s + 1}(dq).
\end{equation}

Let $ E = \mathbb R $ and $A_w  =( - \infty, \, w] $, where
$ w \in \mathbb R $. Denote  $ F(s,x,y,z,t,w) =
P(s,x,y,z,t,A_w ) $. It is clear that $ F(s,x,y,z,t,w) $ as the function of
$ w $ is monotone, right-continuous and $ F(s,x,y,z,t, - \infty ) = 0$, $ F(s,x,y,z,t, +\infty) = 1.$

The condition (V) of Definition \ref{d2} for the function $ F(s,x,y,z,t,w) $ has the following form
$$ F(s,x,y,z,t,w) = \int\limits_E {\int\limits_E {\int\limits_E
{dF(s,x,y,z,\tau ,u)F(\tau ,u,\vartheta ,q,t,w)m_\tau  (d\vartheta
)} } } m_\tau  (dq).$$

If the function $ F(s,x,y,z,t,w) $ is absolutely continuous with respect to variable $w$ then
$$ F(s,x,y,z,t,w) = \int\limits_{ - \infty }^w {f(s,x,y,z,t,u)du},$$
where $ f(s,x,y,z,t,w) $ is a non-negative function and measurable with respect to variables $ x,y,z,w $, moreover it satisfies
the following conditions
$$\int\limits_{- \infty }^\infty  {f(s,x,y,z,t,w)dw}  = 1.$$
\begin{equation}\label{b12}
f(s,x,y,z,t,w) = \int\limits_E {\int\limits_E {\int\limits_E
{f(s,x,y,z,\tau ,u)f(\tau ,u,\vartheta ,q,t,w)m_\tau  (d\vartheta
)m_\tau  (dq)du} } }. \end{equation}

Form  (\ref{b10}) and (\ref{b11}) for $ F(s,x,y,z,t,w)$ we get
$$ \frac{{\partial F(s,x,y,z,t,w)}}{{\partial
t}} = \int\limits_E {\int\limits_E {\int\limits_E {\frac{{\partial
F(s,x,y,z,t - 1,u)}}{{\partial u}}C(t,u,\vartheta ,q,w)dum_{t - 1}
(d\vartheta )m_{t - 1} (dq)} } }, $$

$$ \frac{{\partial F(s,x,y,z,t,w)}}{{\partial s}} = \int\limits_E
{\int\limits_E {\int\limits_E {\frac{{\partial C(s +
1,x,y,z,u)}}{{\partial u}}F(s + 1,u,\vartheta ,q,t,w)dum_{s + 1}
(d\vartheta )m_{s + 1} (dq)} } }. $$

Assume the following limit exists:
$$ a(t,u,\vartheta ,q,w) = \mathop {\lim }\limits_{\Delta
\to 0} \frac{{f(t - 1,u,\vartheta ,q,t + \Delta ,w) -
f(0,u,\vartheta ,q,1,w)}}{\Delta }.$$ Then we obtain the following integro-differential equations:
\begin{equation}\label{b13}\frac{{\partial f(s,x,y,z,t,w)}}{{\partial t}} =
\int\limits_E {\int\limits_E {\int\limits_E {a(t,u,\vartheta
,q,w)f(s,x,y,z,t - 1,u)dum_{t - 1} (d\vartheta )} } } m_{t - 1}
(dq) \end{equation}
\begin{equation}\label{b14}
\frac{{\partial f(s,x,y,z,t,w)}}{{\partial s}} = - \int\limits_E
{\int\limits_E {\int\limits_E {a(s + 1,x,y,z,u)f(s + 1,u,\vartheta
,q,t,w)dum_{s + 1} (d\vartheta )} } } m_{s + 1} (dq).
\end{equation}

\section {Reduction of the integro-differential equations to differential equations }

Under some conditions the integro-differential equations (\ref{b13}) and
(\ref{b14}) can be reduced to differential equations. In this subsection we illustrate this for equation (\ref{b14}).

Let $ t > s + 2 $, then from (\ref{b12}) we get
\begin{equation}\label{b15}
\begin{array}{l}f(s,x,y,z,t,w) - f(s + \Delta ,x,y,z,t,w) = \int\limits_E
{\int\limits_E {\int\limits_E {\left\{ {f(s,x,y,z,s + 1 + \Delta ,u) - } \right.} } }  \\
\left. { - f(s + \Delta ,x,y,z,s + 1 + \Delta ,u)} \right\}f(s + 1 + \Delta ,
u,\vartheta ,q,t,w)m_{s + 1 + \Delta } (d\vartheta )m_{s + 1 + \Delta } (dq)du \\
 \end{array}.
 \end{equation}
We consider function $f$ such that the decomposition of $ f(s + 1 + \Delta ,u,\vartheta ,q,t,w) $ into Taylor's series in a
neighborhood of the point $ (x,y,z)$ has the form:
$$ f(s + 1 +
\Delta ,u,\vartheta ,q,t,w) = f(s + 1 + \Delta ,x,y,z,t,w) +
\frac{{\partial f(s + 1 + \Delta ,x,y,z,t,w)}}{{\partial u}}(u -
x)$$
$$
+ \frac{{\partial f(s + 1 + \Delta ,x,y,z,t,w)}}{{\partial
\vartheta }}(\vartheta  - y) + \frac{{\partial f(s + 1 + \Delta
,x,y,z,t,w)}}{{\partial q}}(q - z)$$
$$
+ \frac{1}{2}\frac{{\partial ^2 f(s + 1 + \Delta
,x,y,z,t,w)}}{{\partial u^2 }}(u - x)^2  +
\frac{1}{2}\frac{{\partial ^2 f(s + 1 + \Delta
,x,y,z,t,w)}}{{\partial \vartheta ^2 }}(\vartheta  - y)^2$$
$$
+ \frac{1}{2}\frac{{\partial ^2 f(s + 1 + \Delta
,x,y,z,t,w)}}{{\partial q^2 }}(q - z)^2  + \frac{{\partial ^2 f(s
+ 1 + \Delta ,x,y,z,t,w)}}{{\partial u\partial \vartheta }}(u -
x)(\vartheta  - y)$$
$$
+ \frac{{\partial ^2 f(s + 1 + \Delta ,x,y,z,t,w)}}{{\partial
u\partial q}}(u - x)(q - z) + \frac{{\partial ^2 f(s + 1 + \Delta
,x,y,z,t,w)}}{{\partial \vartheta \partial q}}(\vartheta  - y)(q -
z)$$
$$
+ \frac{1}{2}\frac{{\partial ^3 f(s + 1 + \Delta
,x,y,z,t,w)}}{{\partial u^2 \partial \vartheta }}(u - x)^2
(\vartheta  - y) + \frac{1}{2}\frac{{\partial ^3 f(s + 1 + \Delta
,x,y,z,t,w)}}{{\partial u^2 \partial q}}(u - x)^2 (q - z)$$
$$
+ \frac{1}{2}\frac{{\partial ^3 f(s + 1 + \Delta
,x,y,z,t,w)}}{{\partial \vartheta ^2 \partial u}}(\vartheta  -
y)^2 (u - x) + \frac{1}{2}\frac{{\partial ^3 f(s + 1 + \Delta
,x,y,z,t,w)}}{{\partial \vartheta ^2 \partial q}}(\vartheta  -
y)^2 (q - z)$$
$$
+ \frac{1}{2}\frac{{\partial ^3 f(s + 1 + \Delta
,x,y,z,t,w)}}{{\partial q^2 \partial u}}(q - z)^2 (u - x) +
\frac{1}{2}\frac{{\partial ^3 f(s + 1 + \Delta
,x,y,z,t,w)}}{{\partial q^2 \partial \vartheta }}(q - z)^2
(\vartheta  - y)$$
$$
+ \frac{1}{6}\frac{{\partial ^3 f(s + 1 + \Delta
,x,y,z,t,w)}}{{\partial q^3 }}(u - x)^3  +
\frac{1}{6}\frac{{\partial ^3 f(s + 1 + \Delta
,x,y,z,t,w)}}{{\partial \vartheta ^2 }}(\vartheta  - y)^3$$
\begin{equation}\label{b16}
+ \frac{1}{6}\frac{{\partial ^3 f(s + 1 + \Delta
,x,y,z,t,w)}}{{\partial q^2 }}(q - z)^3  + \frac{{\partial ^3 f(s
+ 1 + \Delta ,x,y,z,t,w)}}{{\partial u\partial \vartheta \partial
q}}(u - x)(\vartheta  - y)(q - z).
\end{equation}
Substituting (\ref{b16}) into (\ref{b15}) we consider non-zero summands:
 $$
\int\limits_E {\int\limits_E {\int\limits_E {\left\{ {f(s,x,y,z,s + 1 +
\Delta ,u) - f(s + \Delta ,x,y,z,s + 1 + \Delta ,u)} \right\}} } }\times$$
$$\frac{{\partial f(s + 1 + \Delta ,x,y,z,t,w)}}{{\partial x}}(u - x)
m_{s + 1 + \Delta } (d\vartheta )m_{s + 1 + \Delta } (dq)du$$ $$
= \frac{{\partial f(s + 1 + \Delta ,x,y,z,t,w)}}{{\partial
x}}\int\limits_E {\left\{ {f(s,x,y,z,s + 1 + \Delta ,u)} \right. -
} $$ $$ - \left. {f(s + \Delta ,x,y,z,s + 1 + \Delta ,u)}
\right\}(u - x)du.
$$
Denote
$$a(s,x,y,z,\Delta ) = \int\limits_E {\left\{ {f(s,x,y,z,s
+ 1 + \Delta ,u)} \right.}  - \left. {f(s + \Delta ,x,y,z,s + 1 +
\Delta ,u)} \right\}(u - x)du. $$

Now consider the summands with second order of derivations $$
\int\limits_E {\int\limits_E {\int\limits_E {\left\{ {f(s,x,y,z,s
+ 1 + \Delta ,u) - f(s + \Delta ,x,y,z,s + 1 + \Delta ,u)}
\right\}} } }\times$$ $$ \frac{1}{2}\frac{{\partial ^2 f(s + 1 + \Delta
,x,y,z,w)}}{{\partial x^2 }}(u - x)^2 m_{s + 1 + \Delta } (d\vartheta )m_{s + 1 + \Delta }
(dq)du$$ $$ = \frac{1}{2}\frac{{\partial ^2 f(s + 1 + \Delta
,x,y,z,w)}}{{\partial x^2 }}\int\limits_E^{} {\left\{ {f(s,x,y,z,s
+ 1 + \Delta ,u) - } \right.} $$
$$
- \left. {f(s + \Delta ,x,y,z,s + 1 + \Delta ,u)} \right\}(u -
x)^2 du. $$

Denote $$ b^2 (s,x,y,z,\Delta ) = \int\limits_E^{}
{\left\{ {f(s,x,y,z,s + 1 + \Delta ,u) - \left. {f(s + \Delta
,x,y,z,s + 1 + \Delta ,u)} \right\}(u - x)^2 du} \right.} $$
Then
$$ \int\limits_E {\int\limits_E {\int\limits_E {\left\{
{f(s,x,y,z,s + 1 + \Delta ,u) - f(s + \Delta ,x,y,z,s + 1 + \Delta
,u)} \right\}} } }\times$$ $$ \frac{{\partial ^2 f(s + 1 + \Delta
,x,y,z,w)}}{{\partial x\partial y}}(u - x)(\vartheta  - y)m_{s + 1 + \Delta }(d\vartheta )m_{s + 1 + \Delta }(dq)du$$
$$ =\frac{{\partial ^2 f(s + 1 + \Delta ,x,y,z,t,w)}}{{\partial
x\partial y}}\int\limits_E {\left\{ {f(s,x,y,z,s + 1 + \Delta ,u)
- } \right.} $$
$$
- \left. {f(s + \Delta ,x,y,z,s + 1 + \Delta ,u)} \right\}(u -
x)du\int\limits_E {(\vartheta  - y} )m_{s + 1 + \Delta }
(d\vartheta ). $$
Denote
 $$ \int\limits_E^{} {(\vartheta  -y)m_{s + 1 + \Delta } }
(d\vartheta ) = \alpha (s + 1,y,\Delta ).$$ Consequently
$$\int\limits_E \int\limits_E \int\limits_E \left\{ f(s,x,y,z,s
+ 1 + \Delta ,u) - f(s + \Delta ,x,y,z,s + 1 + \Delta ,u)
\right\}\times $$ $$\frac{{\partial ^2 f(s + 1 + \Delta
,x,y,z,t,w)}}{{\partial x\partial z}}(u - x)(q - z)m_{s + 1 + \Delta }(d\vartheta )m_{s + 1 + \Delta}
(dq)du$$ $$ = \frac{{\partial ^2 f(s + 1 + \Delta
,x,y,z,t,w)}}{{\partial x\partial z}}\int\limits_E{\left\{
{f(s,x,y,z,s + 1 + \Delta ,u)} \right.}$$
$$
 - \left. {f(s + \Delta ,x,y,z,s + 1 + \Delta ,u)} \right\}
\int\limits_E {(u - x)du} \int\limits_E {(q - z)m_{s + 1 + \Delta
} } (dq). $$

 $$ \int\limits_E \int\limits_E \int\limits_E
\left\{f(s,x,y,z,s + 1 + \Delta ,u) - f(s + \Delta ,x,y,z,s + 1
+ \Delta ,u)\right\}\times$$ $$\frac{1}{2}\frac{{\partial ^3 f(s + 1 +
\Delta ,x,y,z,t,w)}}{{\partial x^2 \partial y}}(u - x)^2
(\vartheta  - y)m_{s + 1 + \Delta } (d\vartheta )m_{s + 1 +
\Delta } (dq)du$$ $$ = \frac{1}{2}\frac{{\partial ^3 f(s + 1 + \Delta
,x,y,z,t,w)}}{{\partial x^2 \partial y}}\int\limits_E {\left\{
{f(s,x,y,z,s + 1 + \Delta ,u) - } \right.} $$

$$ - \left. {f(s + \Delta ,x,y,z,s + 1 + \Delta ,u)} \right\}(u -
x)^2 du\int\limits_E {(\vartheta  - y)m} _{s + 1 + \Delta }
(d\vartheta )$$ $$ = \frac{1}{2}\frac{{\partial ^3 f(s + 1 + \Delta
,x,y,z,t,w)}}{{\partial x^2 \partial y}}b^2 (s,x,y,z,\Delta )\alpha (s + 1,y,\Delta ).$$

$$ \int\limits_E\int\limits_E \int\limits_E \left\{
{f(s,x,y,z,s + 1 + \Delta ,u) - f(s + \Delta ,x,y,z,s + 1 + \Delta
,u)} \right\}\times $$ $$\frac{1}{2}\frac{{\partial ^3 f(s + 1 + \Delta
,x,y,z,t,w)}}{{\partial x^2 \partial z}}(u - x)^2(q - z)m_{s + 1+\Delta}(d\vartheta )m_{s + 1 + \Delta }
(dq)du$$ $$ = \frac{1}{2}\frac{{\partial ^3 f(s + 1 + \Delta
,x,y,z,t,w)}}{{\partial x^2 \partial z}}\int\limits_E {\left\{
{f(s,x,y,z,s + 1 + \Delta ) - } \right.} $$
$$ - \left. {f(s + \Delta ,x,y,z,s + 1 + \Delta ,u)} \right\}(u -
x)^2 du\int\limits_E {(q - z)m} _{s + 1 + \Delta } (dq)$$ $$ =
\frac{1}{2}\frac{{\partial ^3 f(s + 1 + \Delta
,x,y,z,t,w)}}{{\partial x^2 \partial z}} b^2 (s,x,y,z,\Delta )\alpha (s + 1,z,\Delta ).$$

$$ \int\limits_E \int\limits_E \int\limits_E \left\{
{f(s,x,y,z,s + 1 + \Delta ,u) - f(s + \Delta ,x,y,z,s + 1 + \Delta
,u)} \right\}\times $$ $$ \frac{1}{2}\frac{{\partial ^3 f(s + 1 + \Delta
,x,y,z,t,w)}}{{\partial ^2 y\partial x}}(\vartheta-y)^2
(u - x)m_{s + 1 + \Delta } (d\vartheta )m_{s + 1 + \Delta }
(dq)du$$ $$ = \frac{1}{2}\frac{{\partial ^3 f(s + 1 + \Delta
,x,y,z,t,w)}}{{\partial ^2 y\partial x}}\int\limits_E {\left\{
{f(s,x,y,z,s + 1 + \Delta ,u) - } \right.} $$

$$
- \left. {f(s + \Delta ,x,y,z,s + 1 + \Delta ,u)} \right\}(u -
x)du\int\limits_E {(\vartheta  - y)^2 m} _{s + 1 + \Delta }
(d\vartheta )$$ $$ = \frac{1}{2}a(s,x,y,z,\Delta ) \frac{{\partial ^3 f(s + 1 + \Delta ,x,y,z,t,w)}}{{\partial ^2
y\partial x}} \cdot \alpha _2 (s + 1,y,\Delta ),$$
where $$
\alpha _2 (s + 1,y,\Delta ) = \int\limits_E {(\vartheta  - y)^2
m_{s + 1 + \Delta } (d\vartheta )}.$$

$$ \int\limits_E \int\limits_E \int\limits_E \left\{
{f(s,x,y,z,s + 1 + \Delta ,u) - f(s + \Delta ,x,y,z,s + 1 + \Delta
,u)} \right\}\times$$ $$\frac{1}{2}\frac{{\partial ^3 f(s + 1 + \Delta
,x,y,z,t,w)}}{{\partial z^2 \partial x}}(q - z)^2(u - x)m_{s + 1 + \Delta } (d\vartheta )m_{s + 1 + \Delta }
(dq)du$$ $$ =\frac{1}{2}\frac{{\partial ^3 f(s + 1 + \Delta
,x,y,z,t,w)}}{{\partial z^2 \partial x}}\int\limits_E {\left\{
{f(s,x,y,z,s + 1 + \Delta ,u)-}\right.}$$

$$
- \left. {f(s + \Delta ,x,y,z,s + 1 + \Delta ,u)} \right\}(u -
x)du\int\limits_E {(q - z)^2 m} _{s + 1 + \Delta } (dv)$$ $$ =
\frac{1}{2}\frac{{\partial ^3 f(s + 1 + \Delta
,x,y,z,t,w)}}{{\partial z^2 \partial x}}a(s,x,y,z,\Delta )\alpha _2 (s + 1,z,\Delta ).$$
$$
\int\limits_E \int\limits_E \int\limits_E \left\{ {f(s,x,y,z,s
+ 1 + \Delta ,u) - f(s + \Delta ,x,y,z,s + 1 + \Delta ,u)}
\right\}\times$$ $$\frac{{\partial ^3 f(s + 1 + \Delta
,x,y,z,t,w)}}{{\partial z\partial y\partial z}}(u - x)(\vartheta  - y)(q - z)m_{s + 1 + \Delta } (d\vartheta )m_{s + 1
+ \Delta } (dq)du$$ $$ = \frac{{\partial ^3 f(s + 1 + \Delta
,x,y,z,t,w)}}{{\partial x\partial y\partial z}}\int\limits_E
{\left\{ {f(s,x,y,z,s + 1 + \Delta ,u) - } \right.} $$

$$
- \left. {f(s + \Delta ,x,y,z,s + 1 + \Delta ,u)} \right\}(u -
x)du\int\limits_E {(\vartheta  - y)m_{s + 1 + \Delta } }
(d\vartheta ) \cdot \int\limits_E {(q - z)m_{s + 1 + \Delta } }
(dq) = $$

$$
= \frac{{\partial ^3 f(s + 1 + \Delta ,x,y,z,t,w)}}{{\partial
x\partial y\partial z}}a(s,x,y,z,\Delta )\alpha (s + 1,y,\Delta
)\alpha (s + 1,z,\Delta ). $$
Since other summands are equal to zero, we get
 $$ f(s,x,y,z,t,w) - f(s + \Delta ,x,y,z,t,w) =
a(s,x,y,z,\Delta )\frac{{\partial f(s + 1 + \Delta
,x,y,z,t,w)}}{{\partial x}}$$
$$
+ b^2 (s,x,y,z,\Delta ) \cdot \frac{1}{2}\frac{{\partial ^2 f(s +
1 + \Delta ,x,y,z,t,w)}}{{\partial x^2 }} +
\frac{1}{2}a(s,x,y,z,\Delta )\alpha (s + 1,y,\Delta ) \times $$
$$ \frac{{\partial ^2 f(s + 1 + \Delta ,x,y,z,t,w)}}{{\partial
x\partial y}} + a(s,x,y,z,\Delta )\frac{{\partial ^2 f(s + 1 +
\Delta ,x,y,z,t,w)}}{{\partial x\partial z}}\alpha (s + 1,z,\Delta
)$$
$$
+ \frac{1}{2}\frac{{\partial ^3 f(s + 1 + \Delta
,x,y,z,t,w)}}{{\partial x^2 \partial y}}b^2 (s,x,y,z,\Delta
)\alpha (s + 1,y,\Delta ) + \frac{1}{2}\frac{{\partial ^3 f(s + 1
+ \Delta ,x,y,z,t,w)}}{{\partial x^2 \partial z}} \times $$
$$ b^2 (s,x,y,z,\Delta )\alpha (s + 1,z,\Delta ) +
\frac{1}{2}a(s,x,y,z,\Delta )\frac{{\partial ^3 f(s + 1 + \Delta
,x,y,z,t,w)}}{{\partial y^2 \partial x}} \times $$
$$
\alpha _2 (s + 1,y,\Delta )  + \frac{1}{2}\frac{{\partial
^3 f(s + 1 + \Delta ,x,y,z,t,w)}}{{\partial z^2 \partial
x}}a(s,x,y,z,\Delta )\alpha _2 (s + 1,z,\Delta ) +
a(s,x,y,z,\Delta ) \times $$
\begin{equation}\label{b17}
\times \alpha (s + 1,y,\Delta )\alpha (s + 1,z,\Delta
)\frac{{\partial ^3 f(s + 1 + \Delta ,x,y,z,t,w)}}{{\partial
x\partial y\partial z}}.
\end{equation}
Dividing both sides of the equality (\ref{b17}) by $\Delta $ and passing to
the limit (assuming the limits exist)
as $\Delta  \to 0$, and denoting  $$
A(s,x,y,z) = \mathop {\lim }\limits_{\Delta  \to 0}
\frac{{a(s,x,y,z,\Delta )}}{\Delta }, \   \ B^2 (s,x,y,z) = \mathop
{\lim }\limits_{\Delta  \to 0} \frac{{b^2 (s,x,y,z,\Delta
)}}{{2\Delta }} , $$ $$  D(s + 1,y) = \mathop {\lim
}\limits_{\Delta \to 0} \alpha (s + 1,y,\Delta )  , \   \ D_2 (s +
1,y) = \mathop {\lim }\limits_{\Delta  \to 0} \frac{{\alpha _2 (s
+ 1,y,\Delta )}}{{2\Delta }}$$
from (\ref{b17}) we get $$ \frac{{\partial f(s,x,y,z,t,w)}}{{\partial
s}} =  - A(s,x,y,z)\frac{{\partial f(s + 1,x,y,z,t,w)}}{{\partial
x}} - B^2 (s,x,y,z)\frac{{\partial ^2 f(s +
1,x,y,z,t,w)}}{{\partial x^2 }}$$
$$
- \frac{1}{2}A(s,x,y,z)D(s + 1,y)\frac{{\partial ^2 f(s +
1,x,y,z,t,w)}}{{\partial x\partial y}}
- A(s,x,y,z)D(s + 1,z)\frac{{\partial ^2 f(s +
1,x,y,z,t,w)}}{{\partial x\partial z}}$$ $$ - B^2 (s,x,y,z)D(s +
1,y)\frac{{\partial ^3 f(s + 1,x,y,z,t,w)}}{{\partial x^2 \partial
y}}$$
$$
- B^2 (s,x,y,z)D(s + 1,z)\frac{{\partial ^3 f(s +
1,x,y,z,t,w)}}{{\partial x^2 \partial z}} -
\frac{1}{2}A(s,x,y,z)D_2 (s + 1,y)\frac{{\partial ^3 f(s +
1,x,y,z,t,w)}}{{\partial y^2 \partial x}}$$
$$
- \frac{1}{2}A(s,x,y,z)D_2 (s + 1,z)\frac{{\partial ^3 f(s +
1,x,y,z,t,w)}}{{\partial z^2 \partial x}}$$ \begin{equation}\label{b18} - A(s,x,y,z)D(s +
1,y)D(s + 1,z)\frac{{\partial ^3 f(s + 1,x,y,z,t,w)}}{{\partial x\partial y\partial z}}.
\end{equation}

\begin{remark} We note that the equation (\ref{b18}) is known as differential equation
with advanced argument. Similarly one can reduce the equation
(\ref{b13}) to a differential equation which is known as a
differential equation with delay argument. For theory of such kind
of equations we refer to \cite{EN} and \cite{Mi}.
\end{remark}
\begin{example} Let $E=\mathbb R$. If $m_\tau(du)=r_\tau(u)du$ then from (\ref{b12}) we get
\begin{equation}\label{b12m}
f(s,x,y,z,t,w) = \int\limits_{- \infty }^\infty {\int\limits_{- \infty }^\infty {\int\limits_{- \infty }^\infty
{f(s,x,y,z,\tau ,u)f(\tau ,u,\vartheta ,q,t,w)r_\tau  (\vartheta
)r_\tau(q)dud\vartheta dq} } }. \end{equation}
Let $a(s,t)$ and $b(t)$ be strictly positive functions. Consider
\begin{equation}\label{f}
f(s,x,y,z,t,w)=(\pi a(s,t))^{-1/2}\exp\left[-{(w-x-y-z)^2\over a(s,t)}\right],\end{equation}
$$r_t(u)=(\pi b(t))^{-1/2}\exp\left[-{u^2\over b(t)}\right].$$
For any $A, B, C\in \mathbb R$, with $A>0$ it is known\footnote{see https://en.wikipedia.org/wiki/Gaussian$_-$integral} that
$$\int\limits_{- \infty }^\infty\exp(-A x^2+Bx+C)dx=\sqrt{{\pi\over A}}\exp\left({B^2\over 4A}+C\right).$$
Using this formula (three times) one can see that the function (\ref{f}) satisfies (\ref{b12m}) iff functions $a(s,t)>0$ and $b(t)>0$ satisfy
the following equation
\begin{equation}\label{ab}
a(s,t)=2b(\tau)+a(s,\tau)+a(\tau,t), \ \ \mbox{for any}\ \ s <\tau <t \ \ \mbox{with} \ \  \tau-s\ge 1, \ \ t -\tau \ge 1.
\end{equation}
Thus we obtain a family of CSP generated by functions (\ref{f}) with $a(s,t)$, $b(t)$ satisfying (\ref{ab}).
For example, take $\epsilon\in (0,1)$, $a(s,t)=t-s-\epsilon$ and $b(t)=\epsilon/2$. These functions satisfy (\ref{ab}) and the corresponding
$f$ and $r$ are defined by
$$f(s,x,y,z,t,w)={\exp\left[-{(w-x-y-z)^2\over t-s-\epsilon}\right]\over \sqrt{(t-s-\epsilon)\pi}},$$
$$r_t(u)={\exp\left[-{2u^2\over \epsilon}\right]\over \sqrt{\epsilon\pi/2}}$$
generate a CSP.
\end{example}
{\bf Acknowledgments.} We thank
the referee for careful reading of the manuscript and for useful suggestions.


\section*{References}
\begin{thebibliography}{9}

\bibitem{A} Avise JC 2014 \textit{The Hope, Hype, and Reality of Genetic Engineering} ( New York: Oxford Univ. Press).

\bibitem{B} Bernshtein S N 1924 \textit{Uch. Zap. Nauchno-Issled.
kaf. Ukr. Otd. Mat.} \textbf{1}  83--115.

\bibitem{D} Dohtani A 1992 \textit{SIAM J.Appl. Math.}
\textbf{52} 1707-–1721.

\bibitem{EN} Elgolz L E and Norkin S B 1971 \textit{ Introduction to the theory of differential
equations with a deviating agument} (Moscow: Nauka)

\bibitem{FG} Fisher M E and Goh B S 1977 \textit{J. Math. Biol.}
\textbf{4} 265–-274

\bibitem{GaR} Ganikhodjaev N N and Rozikov U A 2006 {\it Regul. Chaotic
Dyn.} {\bf 11} 467--473

\bibitem{GN} Ganikhodjaev N N 2000 \textit{Doklady Math.} \textbf{61} 321--323

\bibitem{G} Ganikhodjaev N N 1991 {\it J. Theor. Prob.} {\bf 4} 639-–653

\bibitem{GR1} Ganikhodzhaev R N 1993 \textit{Russian Acad. Sci. Sbornik
Math.} \textbf{76}  489--506

\bibitem{GMR} Ganikhodzhaev R N,  Mukhamedov F M and Rozikov U A
2011 \textit{Inf. Dim. Anal. Quant. Prob. Rel. Fields.}
\textbf{14} 279-335

\bibitem{ES} Erson R C and  Stewart J 1971 \textit{Hum. Hered}.  \textbf{21} 523--542

\bibitem{J} Jenks R 1969 \textit{J. Differ. Equations}.  \textbf{5} 497--514

\bibitem{Kh} Khamraev A Yu 2009 \textit{Uzbek. Mat. Zh.} {\bf 3} 65--71

\bibitem{Kh1} Khamraev A Yu 2004 \textit{Uzbek. Mat. Zh.}   {\bf 2} 79--84.

\bibitem{ly} Lyubich Yu  I 1992 \textit{Mathematical structures in population genetics}
(Berlin: Springer)

\bibitem{mal} Malevanets A and  Kapral R 2000 \textit{J. Chem. Phys.} {\bf 112} 7260--7269

\bibitem{M}  Mamurov B J 2010 \textit{Uzbek. Mat. Zh}. {\bf 4} 113--117

\bibitem{Mi} Mishkis A D 1972 \textit{Linear differensial equations with delay
of argument} (Moscow: Nauka) (Russian)

\bibitem{MG} Mukhamedov F and Ganikhodjaev N 2015 \textit{Lect. Notes Math.} {\bf 2133} (Berlin: Springer)

\bibitem{MS2015}  Mukhamedov F and Supar  N A 2015 \textit{Bull. Malay. Math. Sci.
Soc.} {\bf 38} 1281--1296

\bibitem{MS2013}  Mukhamedov F, Supar N A and Pah Ch H 2013
 \textit{J. Phys.: Conf. Ser.} {\bf 435}  012013.

 \bibitem{RK1} Rozikov U A and Khamraev A Yu 2004 \textit{Ukraine Math. Jour.} \textbf{56} 1699--1711

\bibitem{RK} Rozikov U A and Khamraev A Yu 2014  \textit{ Nonlinear Dyn.
Syst. Theory.} \textbf{14} 92--100

\bibitem{S1} Sarymsakov T A and  Ganikhodzhaev N N 1989 \textit{ Soviet
Math. Dokl.} {\bf 39} 369-–373

\bibitem{S2} Sarymsakov T A and  Ganikhodzhaev N N 1991 {\it Soviet Math. Dokl.} {\bf
43}  279–-283

\bibitem{S3} Sarymsakov T A and  Ganikhodzhaev N N 1990 \textit{ J.
Theor. Prob.} {\bf 3} 51–-70


\bibitem{UR} Udwadia F E and Raju N 1998 \textit{Physica D} \textbf{111} 16-–26
\end{thebibliography}
\end{document}